 \documentclass[12pt]{article} 
\usepackage[letterpaper,margin=1in]{geometry}
\usepackage{amssymb}
\usepackage{amsmath}
\usepackage{graphicx} 
\usepackage{wrapfig} 

\def\R{{\mathbb R}}
\def\pf{\emph{Proof. }}
\def\qed{$\hfill \blacksquare$}

\def\con{{\rm cone\,}}
\def\cone{{\rm cone\,}}
\def\co{{\rm conv\,}}
\def\aff{{\rm aff\,}}
\def\bd{{\rm bd\,}}
\def\rbd{{\rm rbd\,}}
\def\amin{{\rm argmin\,}}

\def\iff{if and only if }

\def\la{\langle}
\def\ra{\rangle}

\def\alp{\alpha}
\def\bet{\beta}
\def\lam{\lambda}

\def\cone{{\rm cone}\,}
\def\conv{{\rm conv}\,}
\def\oconv{{\overline{\rm conv}\,}}
\def\ext{{\rm ext}\,}
\def\ri{{\rm ri}\,}
\def\be{\begin{equation}}
\def\ee{\end{equation}}

\newtheorem{thm}{Theorem}[section]
\newtheorem{lem}[thm]{Lemma} 

\newtheorem{prop}[thm]{Proposition}
\newtheorem{defn}[thm]{Definition}
\newtheorem{ex}[thm]{Example}
\newtheorem{rem}[thm]{Remark}

\newtheorem{fact}[thm]{Fact}

\numberwithin{equation}{section}
\numberwithin{figure}{section}
\numberwithin{table}{section}

\title{Visible Points in Convex Sets and Best
  Approximation}

\author{Frank Deutsch\thanks{Department of Mathematics, The Pennsylvania State
  University, University Park, PA 16802. \newline Email: \texttt{deutsch@math.psu.edu}}
\and 
Hein Hundal\thanks{The Pennsylvania State University,  146 Cedar Ridge Drive,
 Port Matilda, PA 16870 \newline Email: \texttt{hundalhh@yahoo.com }}
\and Ludmil Zikatanov
  \thanks{Department of Mathematics, The Pennsylvania State
    University, University Park, PA 16802. \newline Email:
    \texttt{ludmil@psu.edu}. The third author was supported in part by
    NSF DMS-0810982 and {DoE DE-SC0006903}.}
} 

\begin{document}
\maketitle 
              
\begin{abstract}
  The concept of a \emph{visible point} of a convex set relative to a
  given point is introduced.  A number of basic properties of such
  visible point sets is developed. In particular, it is shown that
  this concept is useful in the study of best approximation, and it
  also seems to have potential value in the study of robotics.

\bigskip
\noindent
2010 Mathematics Subject Classification: 41A65, 52A27.

\noindent
Keywords and phrases: best approximation from convex sets, visible points in convex sets.
\end{abstract}

\section{Introduction}\label{S: intro}

Unless explicitly stated otherwise, throughout this paper $X$ will
always denote a (real) normed linear space, and $C$ a closed convex
set in $X$.  For any two distinct points $x, v$ in $X$, we define
interval notation analogous to that on the real line by
\begin{eqnarray*}
[x, v]:&=& \{\lam x+(1-\lam)v \mid 0\le \lam \le 1 \}, \\
\left[x,v\right[:&=&\{\lam x+(1-\lam)v \mid 0 <\lam \le 1\}, \\
\left]x, v\right]:&=&\{\lam x+(1-\lam)v \mid  0\le \lam < 1\}=\left[ v, x\right[ , and \\
\left]x,v\right[:&=&\{\lam x+(1-\lam)v \mid 0< \lam < 1\}. \\
\end{eqnarray*}
In other words, $[x,v]$ is just the closed line segment joining  $x$ and $v$,  $[x,v[$ is the same line segment but excluding the end point $v$, and $]x,v[$ is the line segment $[x, v]$ with both end points $x$ and $v$ excluded.

\begin{defn}\label{visible point} Let $x\in X$. A point $v\in C$ is said to be \textbf{visible} to $x$ with respect to  $C$ \iff $[x,v]\cap C=\{v\}$ or, equivalently, $[x,v[\, \cap\, C=\emptyset$. The set of all visible points  to $x$ with respect to $C$ is denoted by $V_C(x)$. 
\end{defn}
Thus 
\be\label{eq  1: vis pt}
V_C(x)=\{v \in C\mid [x,v]\cap C=\{v\}\}=\{v\in C\mid [x,v[\,\cap\, C=\emptyset\}.
\ee

Geometrically, one can regard $V_C(x)$ as the ``light'' that would be cast on the set $C$ if there were a light source at the point $x$ emanating in all directions. Alternatively, one can regard  the set $C$ as an ``obstacle'' in $X$, a ``robot''  is located at a point $x\in X$, and the directions determined by the intervals $[x, v]$, where $v\in V_C(x)$, as directions to be avoided by the robot so as not to collide with the obstacle $C$.

In this paper we begin a study of  visible sets.  In Section \ref{S: visibility}, we will give some characterizations of visible sets (see Lemmas \ref{char of vis pts} and \ref{char of vis 2}, and Theorem \ref{strong sep char} below). We show that the visible set mapping $V_C$ satisfies a translation property just like the well-known metric projection $P_C$ (see Lemma \ref{translation} below).  Recall that the generally set-valued \emph{metric projection} (or nearest point mapping) $P_C$ is defined on $X$ by
\[
P_C(x):=\{y\in C \mid \| x-y\|=\inf_{c\in C} \|x-c\| \}.
\]
Those closed convex sets $C$ such that the set of visible points to each point not in $C$ is the whole set $C$ are precisely the affine sets (Theorem \ref{char of affine sets}).  In Section \ref{S: V and BA} we study the connection between visible points and best approximations. Finally, in Section \ref{S: ba from polytope} we consider characterizing best approximations to a point in a Hilbert space from a polytope, i.e., the convex hull of a  finite set of points.

\section{Visibility from Convex Sets}\label{S: visibility}

The first obvious consequence of the definition of visibility is the following.

\begin{lem}\label{singleton}Let $C$ be a closed convex set in $X$. If $x\in C$, then $V_C(x)=\{x\}$.
\end{lem}
This lemma shows that the most interesting case is when $x\in X\setminus C$ and the main results to follow actually require this condition as part of  their hypotheses.
Indeed,  when $x\notin C$, there are additional useful criteria that  characterize visible points. For any set $C$, let $\bd C$ denote the \emph{boundary} of $C$.

Unlike the metric projection, the visibility operator is never empty-valued.  

\begin{lem}\label{nonempty}Let $C$ be a closed convex set in $X$. Then \begin{enumerate} 
\item[{\rm (1)}] $V_C(x)\ne \emptyset$ for each $x\in X$, and
\item[{\rm (2)}] $V_C(x) \subset \bd C$ for each $x\in X \setminus C$. 
\end{enumerate}
\end{lem}

\pf  (1) Let $x\in X$.  By Lemma \ref{singleton} we may assume that $x\notin C$.  Fix any $y\in C$.  Then the interval $[x,y]$ contains points in $C$ (e.g., $y$) and points not in $C$ (e.g., $x$).  Let
\[
\lam_0:=\sup\{ \lam \in [0,1]\mid \lam x  +(1-\lam) y \in C\}.
\]
Since $C$ is closed, it follows that $v_0:=\lam_0 x+(1-\lam_0)y \in C$. Hence  $\lam_0 < 1$, and $[x,v_0]\cap C=\{v_0\}$. That is, $v_0\in V_C(x)$. 

(2) Fix any $x \in X \setminus C$. To show that $v\in \bd C$ for each $v \in V_C(x)$.  If not, then there exists some $v\in V_C(x)$ such that $v\in C\setminus \bd C$. Hence $v$ is in the \emph{interior} of $C$. Thus there must be a subinterval $[v_0, v]$ of the interval $[x, v]$ which lies in $C$.  Hence $[x,v] \cap C\ne \{v\}$, a contradiction to $v\in V_C(x)$.  \qed

\begin{lem}\label{char of vis pts} {\rm (Characterization of visible points)} Let $C$ be a closed convex set in $X$, $x\in X\setminus C$,  and $v\in C$.  Then the following statements are equivalent:
\begin{enumerate}
\item[{\rm(1)}] $v$ is visible to $x$ with respect to  $C$.
\item[{\rm (2)}] $\lam x+(1-\lam)v \notin C$ for each $0< \lam \le 1$.
\item[{\rm (3)}] $\max \{ \lam \in [0,1] \mid \; \lam x+(1-\lam)v \in C \} =0$.
\end{enumerate}
\end{lem}

\pf $(1) \Rightarrow (2)$.  If (1) holds, then $[x, v[ \, \cap\, C=\emptyset$. Since $[x,v[ \,=\{\lam x+(1-\lam)v \mid 0< \lam \le 1\}$, (2) follows. 

$(2) \Rightarrow (3)$. 
Since $v\in C$, (3) is an obvious consequence of (2).

$(3)  \Rightarrow (1)$. If (3) holds, then $[x, v[\,\cap\, C =\emptyset$.  That is, $v \in V_C(x)$.  \qed

\medskip
Simple examples in the Euclidean plane (e.g., a box) show that although $C$ is convex, $V_C(x)$ is not convex in general. These simple examples also might seem to indicate that $V_C(x)$ is always closed.  However, the following example in 3 dimensions shows that this is false in general.

Consider the subset of Euclidean 3-space $\ell_2(3)$ defined by 
\be\label{not closed}
C:=(1,0,0)+\cone\{(1,\alp,\bet) \mid \alp^2+(\bet-1)^2\le 1\}.
\ee

\begin{ex} The set $C$ defined by {\rm(\ref{not closed})} is  a closed convex subset of $\ell_2(3)$ such that $0\notin C$ and $V_C(0)$ is not closed.
\end{ex}

\begin{figure}[!htb]
  \begin{center}
    \includegraphics[width=0.48\textwidth]{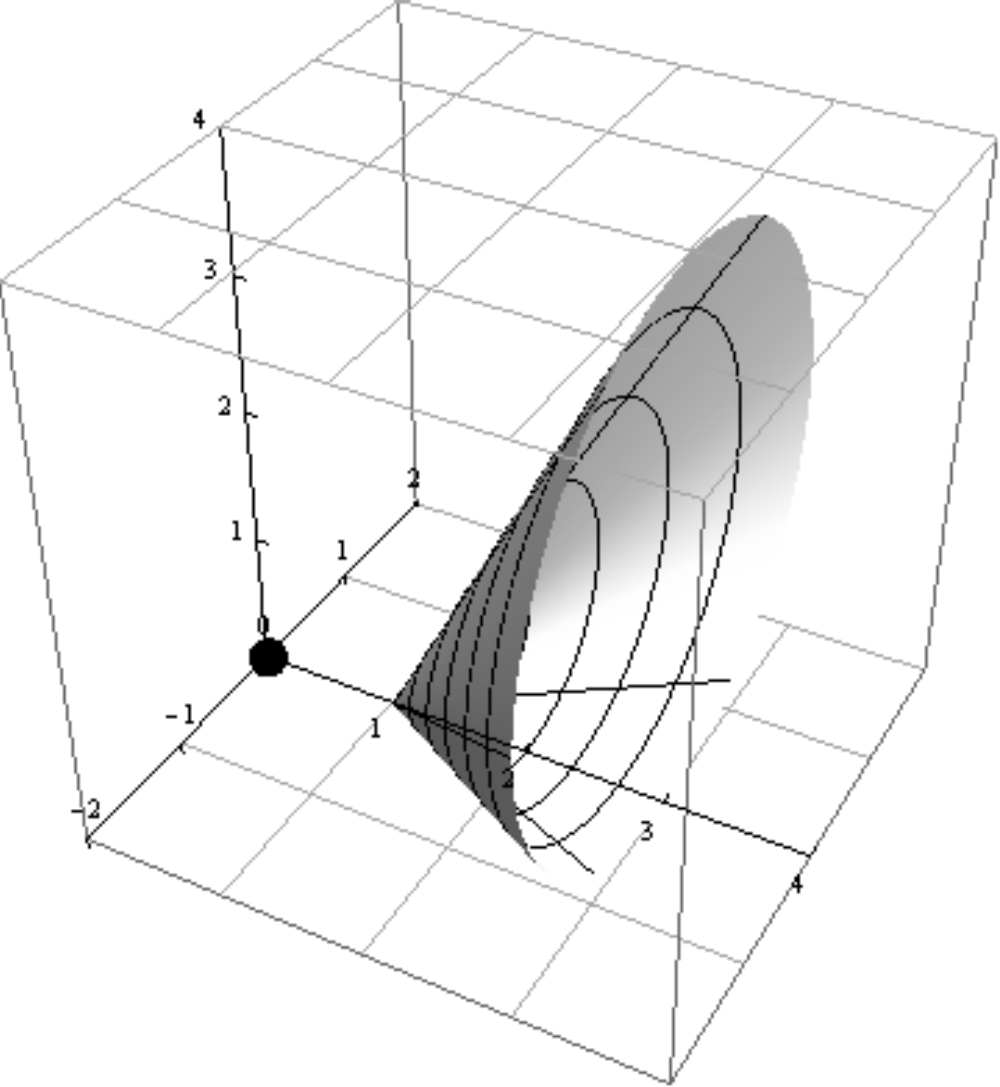}
  \end{center}
  \caption{The set $C$ from Example~2.4.\label{fig24}}
\end{figure}

\pf The result is geometrically obvious (see Figure~\ref{fig24}) 
by observing that the points $(2, \sin t, 1+\cos t)$ are in $V_C(0)$ for
each $0< t< \pi$, but that the limit point $(2,0,0)$ (as $t\to \pi$)
is not. However, the formal proof  
of this fact is a bit
lengthy. Clearly, $0\notin C$ since the first component of any element
of $C$ is at least 1. We first verify the following claim.

\textbf{Claim.} \textsl{The points $v(t):=(2,\, \sin t, \, 1+\cos t)$ are  in $V_C(0)$ for each $0< t < \pi$.}

Using the classical trig identity $\sin^2t+\cos^2t=1$, it is clear that $v(t)\in C$ for each $0< t< \pi$. To complete the proof of the claim, it is enough to show that $[0, v(t)[ \, \cap\, C=\emptyset$ for each $0<t<\pi$. By way of contradiction, suppose the claim is false. Then there exists $0< t_0<\pi$ such that $[0, v(t_0)[ \, \cap\, C\ne \emptyset$. Since $0\notin C$, it follows that there exists $0< \lam<1$ such that $\lam v(t_0) \in C$. That is,
\begin{eqnarray*}
\lam(2, \,\sin t_0, 1+\cos t_0)\in C&=& (1,0,0) +\cone\{(1, \alp, \bet) \mid \alp^2+(\bet-1)^2 \le 1\} \\
&=&(1,0,0)+ \cup_{\rho \ge 0}\rho\{(1, \alp, \bet) \mid \alp^2+(\bet-1)^2 \le 1\}.
\end{eqnarray*}

Since $\lam \sin t_0 \ne 0$, it follows that  for some $\rho>0$, 
\be\label{eq 1; not closed}
\lam (2,\, \sin t_0, 1+\cos t_0)=(1, 0,0)+\rho (1,\alp, \bet)
\ee
  for some $\alp$ and $\bet$ such that 
\be\label{eq 2; not closed}
\alp^2+(\bet -1)^2 \le 1.
\ee
By equating the corresponding components in (\ref{eq 1; not closed}), we obtain
\be\label{eq 3; not closed}
2\lam=1+\rho
\ee
\be\label{eq 4; not  closed}
\lam \sin t_0=\rho \alp
\ee
\be\label{eq 5; not closed}
\lam(\cos t_0+1)=\rho \bet
\ee

From (\ref{eq 3; not closed}) is deduced that $\rho =2\lam -1<2-1=1$ and hence that 
\be\label{eq 6; not closed}
0<\rho<1.
\ee
Also, from (\ref{eq 4; not closed}) and (\ref{eq 5; not closed}) we deduce that $\alp=\mu \sin t_0$ and $\bet=\mu (1+\cos t_0)$, where $\mu:=\lam/\rho$.  Substituting these values for $\alp$ and $\bet$ into (\ref{eq 2; not closed}), we deduce after some algebra that $1 \ge 2\mu^2(1+\cos t_0)-2\mu(1+\cos t_0) +1$.
Subtracting 1 from both sides of this inequality and then dividing both sides of the resulting inequality by the positive number $2\mu(1+\cos t_0)$, we obtain $\mu \le 1$, i.e., $\lam \le \rho$. From (\ref{eq 3; not closed}), it follows  that $\rho \ge 1$, which contradicts (\ref{eq 6; not closed}). This proves the Claim.

It remains to note that the limit point $\lim_{t\to \pi}v(t)=v(\pi)=(2,0,0)$ is \emph{not} in $V_C(0)$.  For this, it enough to note that $[0, v(\pi)[\,\cap\, C\ne \emptyset$.  And for this, it suffices to show that $(3/4)v(\pi)\in C$. But
\[
\hspace{.5in}\frac34v(\pi)=\left(\dfrac64, 0,0\right)=(1,0,0) + \frac12(1,0,0) \in C. \hspace{1.5in} \blacksquare
\]

The following simple fact will be useful to us.  It shows that the visible set mapping $V_C$ satisfies a translation property that is also satisfied by the (generally set-valued) metric projection $P_C$.
 
 \begin{lem}\label{translation} Let $C$ be a closed convex set and $x, y\in X$. Then
 \be\label{eq 1; lem}
 V_C(x)=V_{C+y}(x+y)-y.
 \ee
 \end{lem}
 
 \pf Let $v\in C$. Note that $v\in V_C(x)$ $\Leftrightarrow$  $[x,v[ \,\cap\, C=\emptyset$ $\Leftrightarrow$ $[x+y, v+y[\, \cap\, (C+y)=\emptyset$ $\Leftrightarrow$ $v+y\in V_{C+y}(x+y)$ $\Leftrightarrow$ $v \in V_{C+y}(x+y)-y$.  \qed
 
 \medskip
 It is natural to ask which closed convex sets $C$ have the property that $V_C(x)=C$ for each $x\notin C$.  That is, for which sets is the  whole set  visible to any point outside the set?  The next result shows that this is precisely the class of  affine sets.    Recall that a set $A$ is \emph{affine} if the line through each pair of points in $A$ lies in $A$.  That is, if the line $\aff \{a_1, a_2\}:=\{\alp_1 a_1+\alp_2 a_2\mid \alp_1+\alp_2=1\} \subset A$ for each pair $a_1, a_2 \in A$. Equivalently, $A$ is affine \iff $A=M+a$ for some (unique) linear subspace $M$ (namely, $M=A-A$) and (any) $a\in A$. Finally, the \emph{affine hull} of a set $C$, $\aff(C)$,  is the intersection of all affine sets which contain $C$. As is well-known,
 \be\label{aff hull}
 \aff(C)=\left\{ \sum_{j\in J}\alp_jx_j\biggm| J \mbox{ finite,  }  \sum_{j\in J}\alp_j=1, \, x_j \in C\; \right \}
 \ee
 
 \begin{prop}\label{char of affine sets} Let $C$ be a closed convex set in $X$. Then the following statements are equivalent:
 \begin{enumerate}
 \item[{\rm (1)}] $C$ is affine.
 \item[{\rm (2)}] $V_C(x)=C$ for each $x\in X \setminus C$.
 \end{enumerate}
  \end{prop}
 
 \pf   $(1) \Rightarrow (2)$.  Let us assume first that $C=M$ is actually a subspace, i.e., that  $0\in C$.   Fix  any $x\notin M$. Since $V_M(x) \subset M$, it suffices to show that $M\subset V_M(x)$.  To this end, let $m\in M$. If $m \notin V_M(x)$, then $[x, m[ \, \cap M \ne \emptyset$. Hence there exists $0<\lam <1$ such that $\lam x+(1-\lam)m\in M$. Since $m\in M$, this implies that $\lam x\in M$ and hence $x\in M$, a contradiction. This proves (2) in case $C$ is a subspace.
 
 In general, suppose $C$ is affine. Then $C=M+c$ for some subspace $M$ and $c\in C$. For any $x\in X\setminus C$, we see that $x-c\notin M$ and by the above proof and Lemma \ref{translation} we obtain
 \[
 V_C(x)=V_{M+c}(x)=V_M(x-c)+c=M+c=C.
 \]

 $ (2) \Rightarrow (1)$.  Assume (2) holds.  If $C$ is not affine, then there exist distinct points  $c_1, c_2$ in $C$ such that $\aff \{c_1, c_2\} \not \subset C$.  Since $C$ is closed convex and $\aff\{c_1, c_2\}$ is a line, it follows that either $\aff\{c_1, c_2\}\cap C=[y_1, y_2]$ or $\aff\{c_1, c_2\}\cap C=y_1+\{\rho (y_2-y_1) \mid \rho \ge 0\}$ for some distinct points $y_1, y_2$ in $C$.  In either case, it is easy to verify that $x:=2y_1-y_2 \notin C$. Also, $y_1=\frac12 x+ \frac12 y_2\in [x, y_2[ \, \cap\, C$, which proves that $y_2\notin V_C(x)$ and hence contradicts the hypothesis that $V_C(x)=C$. Thus $C$ must be affine.
  \qed

 \begin{defn}\label{cone gen} Let $C$ be a closed convex subset of $X$. For any point $y\in X$, we define the \textbf{translated cone} $C_y$ of $C$ by
 \[
 C_y:=\cone(C-y) +y.
 \]
 \end{defn}
 
 Some basic facts about the translated cone follow.
 
 \begin{lem}\label{trans cone props} Let $C$ be a closed convex set in $X$.  Then the following statements hold:
 \begin{enumerate}
 \item[{\rm (1)}] $C_y \supset C$ for each $y \in X$.
 \item[{\rm (2)}] The set $\cone(C-y)$, and hence also $C_y$,  is not closed in general.
 \item[{\rm (3)}] If $y\in C$ and the set $\cone (C-y)$ is closed, then $C_y=T_C(y) + y$, where $T_C(y)$ is the \emph{tangent cone} to $C$ at $y$.
 \end{enumerate}
 \end{lem}
 
 \pf  
 (1) $C_y=\con (C-y) +y \supset C-y+y=C$.

 (2) Consider the closed  ball $C$ of radius one in the Euclidean plane centered at the  point $(0,1)$ and let $y$ denote the origin $(0,0)$. Then $C_y$ is the open upper half-plane plus the origin, which is not closed.  
 
 (3) This follows since the definition of the tangent cone to $C$ at the point $y\in C$ is given by $T_C(y)=\overline{\cone} (C-y)$ (see, e.g., \cite[p. 100]{bc;11}).  \qed

\medskip
 One can also characterize the visible points via the translated cone.
 
 \begin{lem}\label{char of vis 2} Let $C$ be a closed convex set in $X$,  $x\in X\setminus C$, and $v\in C$. Then $v\in V_C(x)$ if and only if $x\notin C_v$. Equivalently, $v\notin V_C(x)$ if and only if $x\in C_v$.
 \end{lem}
 
 \pf  If $v\notin V_C(x)$, then $[x,v[ \,\cap\, C\ne \emptyset$. Thus there exists $0< \lam <1$ such that $y:=\lam x+(1-\lam)v\in C$. Hence $x-v=(1/\lam)(y-v)\in \cone (C-v)$ and therefore $x\in C_v$.
 
 Conversely, if $x\in C_v$, then there exist $\rho \ge 0$ and $y\in C$ such that $x=\rho(y-v)+v=\rho y+(1-\rho)v$.  If $\rho \le 1$, then $x$, being a convex combination of two points in $C$, must lie in $C$, a contradiction. It follows that $\rho>1$ and $y=(1/\rho) x+((\rho-1)/\rho) v \in [x, v[\, \cap\, C$. Thus $[x,v[\, \cap\, C \ne \emptyset$, and so $v \notin V_C(x)$ by (\ref{eq 1: vis pt}).  \qed
 
 \medskip
 The following proposition shows that the translated cones of $C$ form the \emph{external} building blocks for $C$.
 
 \begin{prop}\label{C_y builds C} Let $C$ be a closed  convex set in $X$.  Then
\[
\bigcap_{y \in \bd C}C_y=\bigcap_{y\in C} C_y =\bigcap_{y\in X}C_y=C.
\]
\end{prop}

\pf  By Lemma \ref{trans cone props}, $\cap_{y\in X}C_y \supset C$.  Thus to complete the proof, it suffices to show that $\cap_{y\in \bd C}C_y \subset C$. If not, then there exists $x\in \cap_{y\in \bd C}C_y \setminus C$. Thus $x\in C_y\setminus C$ for each $y\in \bd C$. By Lemma \ref{char of vis 2} $y \notin V_C(x)$ for all $y \in \bd C$.  But $V_C(x) \subset \bd C$ by Lemma \ref{nonempty}(2). This shows that $V_C(x)=\emptyset$, which contradicts Lemma \ref{nonempty}(1).  \qed
 
\medskip 
 A somewhat deeper characterization of visible points is available by using the strong separation theorem. Recall that two sets $C_1$ and $C_2$ in the normed linear space $X$ can be  \emph{strongly separated} by a continuous linear functional $x^*\in X^*$ if
  \be\label{s sep}
 \sup_{y\in C_1}x^*(y) <  \inf_{z\in C_2}x^*(z).
 \ee
 One can also interpret strong separation geometrically.  Suppose $C_1$ and $C_2$ are strongly separated by the functional $x^*$ such that (\ref{s sep}) holds.  Let $b$ be any scalar such  that 
\[
\sup_{y\in C_1}x^*(y) \le b \le \inf_{z\in C_2}x^*(z).
\]
Define the hyperplane $H$ and the (open) half-spaces $H^+$ and $H^-$ by
\begin{eqnarray*}
& &H:=\{y \in X \mid x^*(y)=b \},   \quad H^+:=\{y\in X \mid x^*(y)> b\}, \mbox{ and \qquad  }\\
& & H^-: =\{y \in X\mid x^*(y)< b\}.
\end{eqnarray*}
(Note that $H$, $H^-$, and $H^+$ are disjoint sets such that $X=H\cup H^-\cup H^+$.)
Then $H$ is said to \emph{strongly separate} the sets $C_1$ and $C_2$ in the sense that $C_1\subset H\cup H^-$,  $C_2\subset H\cup H^+$, and (at least) one of the sets $C_1$ or $C_2$ is disjoint from $H$.  
 
 \begin{fact}\label{strong sep thm}{\rm(Strong Separation Theorem; see \cite[Theorem V.2.10, p. 417]{ds;58})} Let $C_1$ and $C_2$ be two disjoint closed convex sets in $X$, one of which is compact. Then the sets can be strongly separated by a continuous linear functional.
  \end{fact}
  
  \begin{defn}\label{ext pt} Let $K$ be a convex subset of $X$. A point $e\in K$ is called an \emph{extreme point} of $K$ if $k_1\in K$, $k_2\in K$, $0< \lam < 1$, and $e=\lam k_1+(1-\lam)k_2$ imply that $k_1=k_2=e$.  The set of extreme points of $K$ is denoted by $\ext K$.
  \end{defn}
  
  The following fact is well-known (see, e.g., \cite[pp. 439--440]{ds;58}), and it will be needed in this section and the next.
  
  \begin{fact}\label{krein-milman}{\rm (Krein-Milman)} Let $K$ be a nonempty compact convex subset of $X$. Then:
  \begin{enumerate}
  \item[{\rm (1)}] $K$ has extreme points and $K$ is the closed convex hull of its extreme points: $K=\overline{\conv}(\ext K)$.
  \item[{\rm (2)}] If $x^*\in X^*$,  then $x^*$ attains  its maximum (resp., minimum) value over $K$ at an extreme point of $K$.
  \end{enumerate}
  \end{fact}

 \begin{thm}\label{strong sep char}{\rm(Another characterization of visible points)} Let $C$ be a closed convex subset of $X$, $x\in X\setminus C$, and $v\in C$.  Then the following statements are equivalent:
 \begin{enumerate}
 \item[{\rm (1)}] $v$ is visible to $x$ with respect to $C$; 
 \item[{\rm (2)}] For each point $y\in ]x,v[$, there exists a functional $x^*\in X^*$ that strongly separates $[x, y]$ and $C$, and $x^*(y)=\max_{z\in [x,y]}x^*(z)$;
 \item[{\rm (3)}] For each point $y\in ]x,v[$, there exists a hyperplane $H=H_y$ that contains $y$ and strongly separates $[x, y]$ and $C$.
 
 \end{enumerate}
 \end{thm}
 
 \pf  $(1) \Rightarrow (2)$. Suppose  $v$ is visible to $x$ from $C$.  Then $[x, v[\, \cap\, C=\emptyset$.  In particular, for each $y\in [x,v[$,  $[x,y]\cap C\subset [x,v[\,\cap\, C=\emptyset$.  Thus $[x,y]$ and $C$ are disjoint closed convex sets, and $[x,y]$ is compact.  By Fact \ref{strong sep thm}, there exists $x^*\in X^*$ such that 
 \be\label{eq 1; strict sep char}
 b:=\sup_{z\in [x,y]}x^*(z) < \inf_{c\in C}x^*(c).
 \ee
 
 To verify (2), it remains to show that $x^*(y)=b$. If $x=y$, this is clear. Thus we may assume that $x\ne y$.
 Since $[x,y]$ is compact, the supremum on the left side of (\ref{eq 1; strict sep char}) is attained. Further, this maximum must be attained at an extreme point of $[x,y]$ by Fact \ref{krein-milman}(2). Since $x$ and $y$ are the only two extreme points of $[x, y]$, we must have  $x^*(x)=b$ or $x^*(y)=b$.   
 
 Suppose $x^*(x)=b$. Since $v\in C$, we have $x^*(v)>b$ by (\ref{eq 1; strict sep char}).  Since $y\in ]x,v[$, there exists $0< \lam <1$ such that $y=\lam x+(1-\lam)v$. Then
 \[
 x^*(y)=\lam x^*(x)+(1-\lam)x^*(v) > \lam b+(1-\lam)b=b,
 \]
 which contradicts the definition of $b$.  Thus the condition $x^*(x)=b$ is not possible, and we must have that $x^*(y)=b$, which verifies (2).
 
$(2) \Rightarrow (3)$. Assume (2) holds.  Let $y\in ]x,v[$.  Choose $x^*\in X^*$ as in (2), and define $H:=\{ z\in X \mid x^*(z)=b \}$, where $b=\max_{z\in [x,y]}x^*(z)$. Then $H$ strongly separates $[x,y]$ and $C$,  $x^*(y)=b$,  and so $y\in H$. Thus (3) holds.

$(3) \Rightarrow (1)$.  Suppose (3) holds but (1) fails.  Then $[x,v[ \,\cap\, C\ne\emptyset$. Choose any $y\in ]x,v[ \,\cap\, C$. By (3), there is a hyperplane $H$ that strongly separates $[x,y]$ and $C$ such that $y\in H$. Writing  $H=\{z\in X \mid x^*(z)=b\}$, we see that $[x,y] \subset \{z\in X\mid x^*(z) \le b\}$, $C\subset \{z\in X \mid x^*(z) > b\}$, and $x^*(y)=b$.  But $y\in C$ and hence $x^*(y) >b$, which is a contradiction. \qed

\section{Visibility and Best Approximation}\label{S: V and BA}

In this section we explore the connection between visibility and best approximation. The first such result states that the set of best approximations to $x$ from $C$ is always contained in the set of visible points  to $x$ with respect to $C$. 
 
 \begin{lem}\label{vis vs ba} Let $C$ be a closed convex subset of $X$. Then $P_C(x) \subset V_C(x)$ for each $x\in X$.
\end{lem}
 
 \pf  The result is  trivial if $P_C(x)=\emptyset$. If $x\in C$, then clearly $P_C(x)=\{x\}$ and $V_C(x)=\{x\} $ by Lemma \ref{singleton}. 
 
 Now suppose $x\in X\setminus C$ and let $x_0\in P_C(x)$. Then $x_0\in C$ so $x_0\ne x$. If $[x,x_0[\,\cap\, C\ne \emptyset$, then there exists $0<\lam <1$ such that $x_\lam:=\lam x+(1-\lam)x_0 \in C$. Hence
 \[
 \|x-x_\lam\|=\|(1-\lam)(x-x_0)\|=(1-\lam)\|x-x_0\|< \|x-x_0\|,
 \]
 which is a contradiction to  $x_0$ being a closest point in $C$  to $x$.  This shows that $[x,x_0[\, \cap\, C=\emptyset$ and hence that $x_0\in V_C(x)$.  \qed
 
\medskip   
 Recall that if $X$ is a strictly convex reflexive Banach space, then each closed convex subset $C$ is \emph{Chebyshev} (see, e.g, \cite{sin;70}).  That is, for each $x\in X$, there is a unique best approximation (i.e., nearest point) $P_C(x)$ to $x$ from $C$.  As is  well-known, the most important example of  a strictly convex reflexive Banach space is a Hilbert space. It is convenient to use the following notation. If $S$ is any subset of $X$, then the \emph{convex hull} of $S$ is denoted by $\co{(S)}$ and the closed convex hull of $S$ is denoted by $\overline{\co}(S)$.
 
 Another such relationship between visibility and best approximation is the following.
 
 \begin{lem}\label{vis vs ba 2}  Let $X$ be a strictly convex reflexive Banach space and $C$ a closed convex subset of $X$. Then $C$ is a Chebyshev set and if $x\in X\setminus C$, then
 \be\label{eq 1; vis vs ba}
 P_C(x)=P_{V_C(x)}(x)=P_{\oconv V_C(x)}(x).
 \ee
 \end{lem}
 
 \pf  By Lemma \ref{vis vs ba}, $P_C(x)\in V_C(x)$. Since $V_C(x)\subset \oconv{V_C(x)} \subset C$, it follows that  $P_C(x)\in P_{V_C(x)}(x)$ and $P_C(x)=P_{\oconv{V_C(x)}}(x)$. Thus $P_{V_C(x)}(x)$ is a singleton and (\ref{eq 1; vis vs ba}) holds. \qed

\medskip
While the Krein-Milman theorem (Fact \ref{krein-milman}(1)) shows that the set of extreme points $\ext C$ of a compact convex set $C$ form the \emph{internal} building blocks of $C$,  the next result shows that the sets $C_e$, where $e\in \ext C$, form the \emph{external} building blocks for $C$. It  is a sharpening of Proposition \ref{C_y builds C} in the special case when the closed convex set $C$ is actually compact.

\begin{thm}\label{building blocks} Let $C$ be a compact convex set in $X$. Then
\be\label{eq 1; bb}
C=\bigcap\{C_e\mid e\in \ext C\}=\bigcap\{C_y\mid y\in C\}.
\ee
\end{thm}

\pf  Using Proposition \ref{C_y builds C}, it
suffices to show that  $\cap \{C_e \mid e\in \ext C\} \subset C$. If not, then there exists $x\in \cap\{C_e \mid e\in \ext C\} \setminus C$.  By Fact \ref {strong sep thm}, there exists $x^*\in X^*$ such that 
\be\label{eq 2; bb}
s: =\sup_{c\in C}x^*(c) < x^*(x).
\ee
By compactness of $C$, the supremum of $x^*$ over $C$ is attained, i.e., there exists $c_0\in C$ such that $x^*(c_0)=s$.  As is easily verified, the set
\be\label{eq 3; bb}
\widetilde{C}=C\cap \{y\in X \mid x^*(y)=s\}
\ee
is extremal in $C$, has extreme points (since it is a closed, hence compact, convex subset of $C$), and each extreme point of $\widetilde{C}$ is an extreme point of $C$ (see, e.g., \cite[pp 439--440]{ds;58}).  Choose any extreme point $\tilde{c}$ in $\widetilde{C}$. Then  $\tilde{c}\in \ext C$. Also, $x\in C_{\tilde{c}}=\cone(C-\tilde{c})+\tilde{c}$ implies that $x=\rho (c-\tilde{c}) +\tilde{c}$ for some $\rho>  0$ and  $c\in C$ (see, e.g., \cite[Theorem 4.4(5), p. 45]{deu;01}).  Hence
\begin{eqnarray*}
s&<&x^*(x)=\rho[x^*(c)-x^*(\tilde{c})] +x^*(\tilde{c}) \le x^*(\tilde{c})=s,
\end{eqnarray*}
which is impossible. This contradiction completes the proof.   \qed

 \begin{prop}\label{convex combo 2} Let  $C$ be a closed convex set in $X$, $x\in X\setminus C$,  and let $x_0\in C$ be a proper convex combination of points $e_i$ in $C$. That is, $x_0=\sum_1^k\lam_ie_i$ for some  $\lam_i > 0$ with $\sum_1^k\lam_i=1$.  If $x_0$ is visible to $x$ with respect to $C$, then each $e_i$ is also visible to $x$.
 \end{prop}
 
 \pf If $k=1$ the result is trivial. Assume that $k=2$. (We will reduce the general case to this case.) 
  
 If the result were false, then we may assume without loss of generality that $e_1$ is not visible to $x$. Thus $]x, e_1[\, \cap\, C \ne \emptyset $.  Hence there exists $0< \mu < 1$ such that $x_1:=\mu x+(1-\mu)e_1\in C$. It follows that 
 \be\label{eq 1:cc1}
 e_1=\frac1{1-\mu}x_1-\frac\mu{1-\mu}x.
 \ee
 Next consider, for each $\rho\in [0,1]$, the expression $x(\rho):=\rho x_1+(1-\rho)e_2$.  Clearly, $x(\rho)\in C$ for all such $\rho$ since both $x_1$ and $e_2$ are in $C$ and $C$ is convex. Omitting some simple algebra, we deduce that
\begin{eqnarray*}
 x(\rho)&=&\rho [\mu x+(1-\mu)e_1]+(1-\rho)e_2\\
 &=&\rho\mu x+(1-\rho\mu)x_0+\rho(1-\mu)e_1+(1-\rho)e_2-(1-\rho\mu)x_0\\
 &=&\rho\mu x+(1-\rho\mu)x_0+[\rho(1-\mu+\lam_1\mu)-\lam_1]e_1+[-\rho(1-\mu+\lam_1\mu)+\lam_1]e_2.
 \end{eqnarray*}
 
 In particular, if we choose 
 \be\label{eq 3:cc1}
 \tilde{\rho}:=\frac{\lam_1}{1-\mu+\lam_1\mu}\, ,
 \ee
 it is not hard to check that $0< \tilde{\rho}<1$. Thus $0< \tilde{\rho}\mu<1$ and 
 \be\label{eq 4:cc1}
 x(\tilde{\rho})=\tilde{\rho}\mu x+(1-\tilde{\rho}\mu)x_0\in C.
 \ee
 This proves that $x(\tilde{\rho})\in\, ]x, x_0[\, \cap\, C$, which contradicts the fact that $x_0$ is visible to $x$. 
 
 Finally, consider the case when $k \ge 3$. If the result were false, then without loss of generality, we may assume that $e_1$ fails to be visible to $x$. Write 
 \[
 x_0=\lam_1 e_1+\mu\sum_{i=2}^k\frac{\lam_i}{\mu}e_i,
  \]
 where $\mu:=\sum_2^k\lam_i=1-\lam_1$. Then $0<\mu<1$, $\lam_1=1-\mu$,  and $x_0=(1-\mu)e_1+\mu y$, where $y=\sum_{2}^k\frac{\lam_i}{\mu}e_i\in C$ by convexity.  By the case when $k=2$ that we proved above,  we get that  $e_1$ (as well as $y$) is visible to $x$, which is a  contradiction.  \qed

\begin{rem} Simple examples in the plane (e.g., a triangle) show that the converse to Proposition \ref{convex combo 2} is \emph{false}!  That is, one could have a closed convex set $C$, a point $x\in X\setminus C$, points $e_i \in V_C(x)$ for $i=1,2, \dots, k$, $k\ge 2$, but $x_0=\frac{1}{k }\sum_1^k  e_i \in C$ is not visible to $x$.
\end{rem}

 \begin{thm}\label{bas in fd spaces} Let $C$ be a closed and bounded convex set in an $n$-dimensional normed linear space $X$ such that  $\ext C$ is closed. Then
 \be\label{eq 1; bas}
 C=\left\{\sum_1^{k}\lam_ie_i \biggm| 1\le k \le n+1,\; \lam_i \ge 0, \;\sum_1^{k}\lam_i=1, \;  e_i \in \ext C \right\}.
 \ee
 Further, let  $x\in X\setminus C$. Then  each point in $P_C(x)$ is a proper convex combination of no more than $n+1$ extreme points of $C$ all of which are visible to $x$ with respect to $C$. That is, 
 \be\label{eq 2; bas}
 P_C(x) \subset \left\{\sum_1^{k}\lam_ie_i \biggm| 1\le k \le n+1,\; \lam_i \ge 0, \;\sum_1^{k}\lam_i=1, \; e_i \in (\ext C) \cap V_C(x) \right\}.
 \ee
 \end{thm}
 
 \pf  Since closed and bounded sets in finite-dimensional spaces are compact,  Fact \ref{krein-milman}(1) implies that $C=\overline{\conv }(\ext C)$.  By Caratheodory's theorem (see, e.g., \cite[p. 17]{che;66}), each point in $\conv(\ext C)$ may be expressed as a convex combination of at most $n+1$ points of $\ext C$. That is,
 \be\label{eq 3; bas}
 \conv(\ext C)=\left\{\sum_1^{n+1}\lam_ie_i \biggm| e_i\in \ext C, \lam_i \ge 0, \sum_1^{n+1}\lam_i=1\right\}
 \ee
 We will show that $\conv(\ext C)$ is closed, and hence that $C=\conv(\ext C)$. To this end, let $(x_k)$ be a sequence in $\conv(\ext C)$ such that $x_k \to x$.  It suffices to show that $x\in \conv(\ext C)$. We have that $x_k=\sum_{i=1}^{n+1} \lam_{ki}e_{ki}$ for some $e_{ki}\in \ext C$,  $\lam_{ki}\ge 0$, and $\sum_{i=1}^{n+1}\lam_{ki}=1$ for each $k$.
 
 Since $C$ is compact, and the $\lam_{ki}$ are bounded, it follows that by passing to a subsequence that $e_{ki}\to e_i$ and $\lam_{ki} \to \lam_{i}$ for each $i$.  Since $\ext C$ is closed, $e_i\in \ext C$.  Further, it easy to see that $\lam_i\ge 0$ for each $i$ and $\sum_{i=1}^{n+1}\lam_i=\lim_k\sum_{i=1}^{n+1}\lam_{ki}=\lim_k1=1$. Thus, $x=\sum_1^{n+1}\lam_ie_i \in \conv(\ext C)$ and $\conv(\ext C)$ is closed. This proves (\ref{eq 1; bas}).
 
Now let $x\in X\setminus C$. By the first part, each point of $P_C(x)$ is in $\conv(\ext C)$. By Lemmas \ref{convex combo 2} and \ref{vis vs ba}, (\ref{eq 2; bas}) follows.  \qed

 \section{Best Approximation from a Simplex}\label{S: ba from polytope}

 In this section we investigate the problem of finding best approximations from a \emph{polytope}, i.e.,  the convex  hull of a finite number of points in a  Hilbert space $X$.  Such sets are compact (because they are closed and bounded  in a finite-dimensional subspace). 
  
 Let $E:=\{e_0, e_1, \dots, e_n\}$ be a set of $n+1$  points in $X$ that is \emph{affinely independent}, i.e., $\{e_1-e_0, e_2-e_0, \dots, e_n-e_0\}$ is linearly independent.  This implies that each point in the convex hull  $C=\conv\{e_0, e_1, \dots, e_n\}$ has a unique representation as a convex combination of the points of $E$. In this case, $C$ is also called an \emph{$n$-dimensional simplex} with vertices $e_i$, since the dimension of the affine hull $\aff(C)$ of $C$ is $n$.  Further, the \emph{relative interior} of $C$, that is, the interior of $C$ relative to $\aff(C)$,  is given by
 \be\label{rel int}
 \ri(C):=\biggl\{  \sum_{i=0}^n\lam_ie_i \bigm | \lam_i>0, \; \sum_{i=0}^n\lam_i=1\biggr\}.
 \ee
 It follows that the \emph{relative boundary} of $C$, $\rbd(C):=C\setminus \ri(C)$,  is given by
 \be\label{rel bdry}
 \rbd(C)=\biggl\{ \sum_{i=0}^n\lam_ie_i \bigm |  \lam_i \ge 0, \; \sum_{i=0}^n\lam_i=1, \; \lam_j=0 \mbox{ for at least one $j$} \biggr\}.
 \ee
 (See \cite[p. 44ff]{roc;70} and \cite[p. 7ff]{hol;75} for more detail and proofs about the facts stated in this paragraph.)
 
 We consider sets of affinely independent points, since this case  captures the essence of our constructions and arguments.  Convex hulls of  $n$ affinely dependent points (i.e., finite point sets that are not affinely independent) can be split into the union of a finite number of convex hulls of subsets of affinely independent points. Thus the problem of finding best approximation from the convex hull of an affinely dependent set of points can be reduced to a finite number of  problems analogous to the one that we consider below in detail.

 \medskip
 \textbf{Under the above hypothesis that $C$ is an $n$-dimensional simplex, we wish to compute $P_C(x)$ for any $x\in X$.}
 \medskip 
 
 We  give an explicit formula for $P_C(x)$ in the case when $n=1$; that is, when $C=[e_0, e_1]$ is a line segment. Then, by a recursive argument, we will indicate how to compute $P_C(x)$ when $C$ is an $n$-dimensional simplex for any $n\ge 2$. First we recall that the \emph{truncation function} $[\cdot]_0^1$ is defined on the set of real numbers by
  \[
 [\alp]_0^1=\left \{ \begin{array}{ll} 0 & \mbox{ if $\alp < 0$}\\
                                                            \alp & \mbox{ if $ 0\le \alp \le 1$}\\
                                                            1 & \mbox{ if $\alp > 1$}. 
                                                            \end{array}                       \right.
  \]
 (Note that in the space $X=\R$, $[\alp]_0^1=P_{[0,1]}(\alp)$ for all $\alp\in \R$.)
 \begin{prop}\label{case n=1} Let  $C=\conv\{ e_0, e_1\}=[e_0, e_1]$ be a $1$-dimensional simplex. Then, for each $x\in X$, 
 \be\label{eq 1: ba in line}
 P_C(x)=e_0+\left[\frac{\la x-e_0,e_1-e_0\ra}{\|e_1-e_0\|^2}\right]_0^1(e_1-e_0). 
 \ee
  \end{prop}
 
 \pf  Let $\alp:={\la x- e_0, e_1-e_0\ra}{\|e_1-e_0\|^{-2}}$ and $c_0:=e_0+[\alp]_0^1(e_1-e_0)$. Then $c_0\in C$, and by the well-known characterization of best approximations from convex sets in Hilbert space (see, e.g., \cite[p. 43]{deu;01}) it suffices to show that 
 \be\label{eq 2: ba in line}
 \la x-c_0, y-c_0\ra \le 0   \mbox{\quad for each $y\in C$.}
 \ee
 Let $y\in C$. Then $y=e_0+\lam(e_1-e_0)$ for some $\lam \in [0,1]$. Hence
 \begin{eqnarray*}
 \la x-c_0, y-c_0\ra &= &\la x-e_0-[\alp]_0^1(e_1-e_0), \lam(e_1-e_0)-[\alp]_0^1(e_1-e_0)\ra \\
 &=& (\lam-[\alp]_0^1)[\la x-e_0, e_1-e_0\ra -[\alp]_0^1\|e_1-e_0\|^2] \\
 &=& (\lam -[\alp]_0^1)\|e_1-e_0\|^{2}\left[ \alp-[\alp]_0^1\right] .
  \end{eqnarray*}
  By considering the three possible cases: $\alp<0$, $\alp\in [0,1]$, and $\alp>1$, it is easy to see that the last expression is always $\le 0$.  Hence (\ref{eq 2: ba in line}) is verified. \qed

 \medskip
 Before considering  the cases when $n\ge 2$, let us first consider the problem of computing $P_{A}(x)$ for any $x\in X$, where $A=\aff C$.
 
 \begin{fact}\label{aff ba's} Let $C=\conv\{e_0, e_1, \dots, e_n\}$ be an $n$-dimensional simplex, and let $A=\aff(C)$. For any $x\in X$, we have
 \be\label{aff ba's: eq 1}
 P_A(x)=e_0+\sum_{j=1}^n\alp_j(e_j-e_0), 
 \ee
 where the scalars $\alp_i$ satisfy the ``normal'' equations:
 \be\label{aff ba's: eq 2}
 \sum_{j=1}^n\alp_j \la e_j-e_0, e_i-e_0 \ra=\la x-e_0, e_i-e_0 \ra   \qquad     (i=1,2, \dots, n).
 \ee
 \end{fact}
 
 The proof of this fact can be found e.g., in \cite[p. 418]{bc;11} or \cite[p. 215]{deu;01}.  Moreover, the ``reduction principle'' that was established  in \cite[p. 80]{deu;01} (where it was stated in the particular case of a subspace) can be easily extended to affine sets as follows.
 
 \begin{fact}\label{red prin}{\rm(Reduction Principle)} Let $C$ be a closed convex set in the Hilbert space $X$ and let $A=\overline{\aff}(C)$. Then $P_C=P_C\circ P_A$. That is, for each $x\in X$,
 \[
 P_C(x)=P_C(P_A(x))  \mbox{ \qquad and  \qquad  }   d^2(x, C)=d^2(x, A)+d^2(P_A(x), C).
 \]
 \end{fact}
 
 We are going to use the Reduction Principle  as follows.  We assume that it is straightforward to find the best approximation to any $x$ in the  set $A=\aff C$, where $C$ is an $n$-dimensional simplex (since it involves only solving a linear system of $n$ equations in $n$ unknowns by Fact \ref{aff ba's}).  The Reduction Principle states that (by replacing $x$ with $P_A(x)$ if necessary) we may as well assume that our point $x$ is in $A$ to begin with, and we shall do this in what follows. We will see that the case when $n=2$ can be reduced to the case when $n=1$ (i.e., Proposition \ref{case n=1} above) for which there is an explicit formula.
 
 \begin{prop}\label{case n=2} Let  $C=\conv\{e_0, e_1, e_2\}$ be a $2$-dimensional simplex. Then for each $x\in \aff(C)$, either $x\in C$ in which case $P_C(x)=x$, or $x\notin C$, in which case
 \be\label{eq 1: case n=2} 
 P_C(x)=P_{[e_i, e_{i+1}]}(x) \mbox{\quad for any  $i\in \{0,1,2\} $ that satisfies }
 \ee
 \be\label{eq 2: case n=2}
 \|x-P_{[e_i, e_{i+1}]}(x)\|=\min_j \|x-P_{[e_j, e_{j+1}]}(x)\|.
 \ee
 {\rm(}Here $e_3:=e_0$.{\rm )}
 \end{prop}
 
 \pf If $x\in C$, then obviously $P_C(x)=x$.  Thus we can assume that $x\in \aff(C)\setminus C$. It follows that $P_C(x)$ must lie on $\rbd C=\cup_{i=0}^2[e_i, e_{i+1}]$.  That is, $P_C(x) \in  [e_i, e_{i+1}]$ for some $i=0,1,$ or $2$. 
  
 \textbf{Claim:} $P_C(x)=P_{[e_i, e_{i+1}]}(x)$ for each $i$ such that  $P_C(x)\in [e_i, e_{i+1}]$.
 
 To see this, we observe that since $P_C(x)\in [e_i, e_{i+1}]$, we have
 \[
 \|x-P_C(x)\|=d(x, C) \le d(x, [e_i, e_{i+1}]) \le \|x-P_C(x)\|
 \]
 which implies that $\|x-P_{[e_i, e_{i+1}]}(x)\|=d(x, [e_i, e_{i+1}])=\|x-P_C(x)\|$. By uniqueness of best approximations from convex sets in Hilbert space, the claim is proved.
 
 If  $k$ is any index such that $\|x-P_{[e_k, e_{k+1}]}(x)\|=\min_j\|x-P_{[e_j, e_{j+1}]}(x)\|$, then it is clear that we must have $P_C(x)=P_{[e_k, e_{k+1}]}(x)$.  \qed

\medskip
Now it appears to be straightforward to apply the idea of Proposition \ref{case n=2} to any $n$-dimensional simplex to describe how to determine $P_C(x)$.   

Let $C=\conv\{e_0, e_1, \dots, e_n\}$ be an $n$-dimensional simplex in $X$ and $x\in \aff(C)$.  If $x\in C$, we have $P_C(x)=x$.  Thus we may assume that $x\in \aff(C) \setminus C$.  It follows that $P_C(x) \in \rbd(C)$. From (\ref{rel bdry}) we see
\be\label{rel bdry 2}
\rbd(C)=\biggl\{ \sum_0^n\lam_ie_i \bigm | \lam_i \ge 0, \; \sum_0^n\lam_i=1, \; \lam_j=0 \mbox{  for some $j$} \biggr \}.
\ee
Since every $y\in \rbd C$ is contained in (at least) one of the sets 
\be\label{C_j}
C_j:=\biggl \{ \sum_{i=0}^n \lam_ie_i \bigm | \lam_i\ge 0 \mbox{  for all $i$,  } \lam_j=0,  \, \mbox{ and } \sum_0^n\lam_i=1 \biggr \},
\ee
it follows that 
\be\label{rel bdry 2}
\rbd C=\bigcup_{j=0}^n C_j.
\ee
Further, each $C_j$ is a simplex of  dimension  $n-1$ in $C$, $P_C(x) \in C_j$ for at least one $j$, and for all such $j$, we have  that
\[
\|x-P_{C}(x)\|=d(x, C) \le \|x-P_{C_j}(x)\| = d(x, C_j)\le \|x-P_C(x)\|.
\]
This implies that equality holds throughout these inequalities, and hence by the uniqueness of best approximations, we have $P_{C}(x)=P_{C_j}(x)$.  If $J=\{j \mid \|x-P_{C_j}(x)\|=\min_i\|x-P_{C_i}(x)\| \}$, then clearly  $P_C(x)=P_{C_j}(x)$  for each $j\in J$.

This discussion suggests the following recursive algorithm for computing $P_C(x)$ when $C=\conv\{e_0, e_1, \dots, e_n\}$ is an $n$-dimensional simplex.  Let $C_j$ be the $(n-1)$-dimensional simplices as defined  in (\ref{C_j}). Let $A=\aff C$, $A_j=\aff C_j$ for each $j=0, 1, \dots, n$, $x\in A\setminus C$, and $x_j=P_{C_j}(x_j)$ for all $j$. The algorithm below defines a function $P(n, x, C)$ which takes as input $n$ and $x$ and the set $C$ and returns the best approximation $P_C(x)$.

\medskip
\textbf{Algorithm.} 
\begin{enumerate}
\item[{\rm (1)}] If $n=1$, then find $P(1, x, C)$ by using the formula given in Proposition \ref{case n=1}.
\item[{\rm (2)}] If $n>1$, then compute $x_j=P_{A_j}(x)$ and $P_{C_j}(x_j)=P(n-1, x_j, C_j)$ for $j=0, 1, \dots, n$.
\item[{\rm (3)}] Set $P_C(x)=P_{C_j}(x_j)$ for any $j\in \amin_k\|x_k-P_{C_k}(x_k)\|$.

\end{enumerate}

\bibliographystyle{plain}

\end{document}